\pdfoutput=1
\RequirePackage{ifpdf}
\ifpdf 
\documentclass[pdftex]{sigma}
\else
\documentclass{sigma}
\fi

\begin{document}

\allowdisplaybreaks

\renewcommand{\PaperNumber}{064}

\FirstPageHeading

\ShortArticleName{Dunkl-Type Operators with Projection Terms}

\ArticleName{Dunkl-Type Operators with Projection Terms\\
Associated to Orthogonal Subsystems\\
in Root System}

\Author{Fethi BOUZEFFOUR}

\AuthorNameForHeading{F.~Bouzef\/four}

\Address{Department of Mathematics, King Saudi University, College of Sciences,\\
P.O.~Box 2455 Riyadh 11451, Saudi Arabia}

\Email{\href{mailto:fbouzaffour@ksu.edu.sa}{fbouzaffour@ksu.edu.sa}}

\ArticleDates{Received April 24, 2013, in f\/inal form October 16, 2013; Published online October 23, 2013}

\Abstract{In this paper, we introduce a~new dif\/ferential-dif\/ference operator $T_\xi$ $(\xi \in
\mathbb{R}^N)$ by using projections associated to orthogonal subsystems in root systems.
Similarly to Dunkl theory, we show that these operators commute and we construct an intertwining operator
between $T_\xi$ and the directional derivative $\partial_\xi$.
In the case of one variable, we prove that the Kummer functions are eigenfunctions of this operator.}


\Classification{33C15; 33D52; 35A22}

\section{Introduction}

In a~series of papers~\cite{Dunkl1,Dunkl2,Dunkl3,Dunkl4}, C.F.~Dunkl builds up the framework for a~theory of dif\/ferential-dif\/ference operators and special functions
related to root systems.
Beside them, there are now various further Dunkl-type operators, in particular the trigonometric Dunkl
operators of Heckman~\cite{Heckman1, Heckman2}, Opdam~\cite{Opdam}, Cherednik~\cite{Cherednik}, and the
important $q$-analogues of Macdonald and Cherednik~\cite{Mac}, see also~\cite{Bouzeffour,Koor}.

\looseness=1
The main objective of this paper is to present a~new class of dif\/ferential-dif\/ference ope\-ra\-tors~$T_\xi$,
$\xi \in \mathbb{R}^N$ with the help of orthogonal projections related to orthogonal subsystems in root
systems.
In other words, our operators follow from Dunkl operator after replacing the usual ref\/lections that exist
in the def\/inition of the operator with their corresponding ortho\-go\-nal projections.
Several problems related to the Dunkl theory arise in the setting of our operators, in particular,
commutativity of $\{T_\xi, \, \xi \in \mathbb{R}^N\}$ and the existence of the intertwining operators.

The outline of the content of this paper is as follows.
In Section~\ref{Sect2}, we collect some def\/initions and results related to root systems and Dunkl operators
which will be relevant for the sequel.
In Section~\ref{Sect3}, we introduce new dif\/ferential-dif\/ference operators $T_\xi$ and we prove the
f\/irst main result.
In Section~\ref{Sect4}, we give an explicit formula for the intertwining operator between $T_\xi$ and the
directional derivative.
In Section~\ref{Sect5}, we study the one variable case.
Finally, in Section~\ref{Sect6} we study the cases of orthogonal subsets in root systems of type~$A_{N-1}$ and~$B_N$.

\section{Dunkl operators}\label{Sect2}

Let us begin to recall some results concerning the root systems and Dunkl operators.
A useful reference for this topic is the book by Humphreys~\cite{Humphreys}.
Let $\alpha \in \mathbb{R}^N\backslash \{0\}$, we denote by~$s_{\alpha}$ the ref\/lection onto the
hyperplane orthogonal to $\alpha$; that is,
\begin{gather*}
s_{\alpha}(x)=x-2\frac{\langle x,\alpha\rangle}{|\alpha|^2}\alpha,
\end{gather*}
where $\langle\cdot ,\cdot \rangle $ denotes the Euclidean scalar product on $\mathbb{R}^N$, and $|x|=\sqrt{\langle
x,x\rangle}$.

A root system is a~f\/inite set $R$ of nonzero vectors in $\mathbb R^N$ such that for any $\alpha \in R$
one has
\begin{gather*}
s_{\alpha}(R)=R,
\qquad
\textrm{and}
\qquad
R\cap\mathbb R\alpha=\{\pm\alpha\}.
\end{gather*}
A positive subsystem $R_+$ is any subset of $R$ satisfying $R = R_+\cup \{-R_+\}$.
The Weyl group $W=W(R)$ (or real f\/inite ref\/lection group) generated by the root system $R\subset
\mathbb R^N$ is the subgroup of orthogonal group $O(N)$ generated by $\{s_\alpha,\, \alpha\in R\}$.
A multiplicity function on $R$ is a~complex-valued function $ \kappa: R \rightarrow \mathbb C $ which
is invariant under the Weyl group $W$, i.e.,
\begin{gather*}
\kappa(\alpha)=\kappa(g\alpha),
\qquad
\forall\, \alpha\in R,
\qquad
\forall\, g\in W.
\end{gather*}
Let $\xi\in \mathbb{R}^N$, the Dunkl operator $\mathcal{D}_{\xi}$ associated with the Weyl group $W(R)$ and
the multiplicity function $\kappa$, is the f\/irst order dif\/ferential-dif\/ference operator:
\begin{gather}
(\mathcal{D}_{\xi}f)(x)=\partial_{\xi}f(x)+\sum_{\alpha\in R_+}
\kappa(\alpha)\langle\alpha,\xi\rangle\frac{f(x)-f(s_{\alpha}x)}{\langle x,\alpha\rangle}.
\label{Dunkl}
\end{gather}
Here $\partial_\xi$ is the direction derivative corresponding to~$\xi$ and~$s_\alpha$ is the orthogonal
ref\/lection onto the hyperplane orthogonal to~$\alpha$.

The Dunkl operator $\mathcal{D}_{\xi}$ is a~homogeneous dif\/ferential-dif\/ference operator of degree~$-1$.
By the $W$-invariance of the multiplicity function $\kappa$, we have
\begin{gather*}
g^{-1}\circ\mathcal{D}_{\xi}\circ g=\mathcal{D}_{g\xi},
\qquad
\forall\, g\in W(R),
\qquad
\xi\in\mathbb R^N.
\end{gather*}
The remarkable property of the Dunkl operators is that the family $\{ \mathcal{D}_\xi,\, \xi\in\mathbb R^N\}$
generates a~commutative algebra of linear operators on the $\mathbb C$-algebra of polynomial functions.

\section{Operators of Dunkl-type}\label{Sect3}

Let $R$ be a~root system.
A subset $R'$ of $R$ is called a~subsystem of $R$ if it satisf\/ies the following conditions:
\begin{itemize}\itemsep=0pt
\item[i)] If $\alpha \in R'$, then $-\alpha \in R'$;
\item[ii)] If $\alpha,\beta \in R'$ and $\alpha+\beta \in R$, then $\alpha+\beta \in R'.$
\end{itemize}
A subsystem $R'$ of a~root system $R$ in $\mathbb{R}^N$ consisting of pairwise orthogonal roots is called
orthogonal subsystem.
In this case the related Weyl group $W(R')$ is a~subgroup of $\mathbb{Z}_2^N$.
For a~vector $\alpha\in\mathbb R^N\setminus\{0\}$, we write
\begin{gather*}
\tau_{\alpha}(x)=x-\frac{\langle x,\alpha\rangle}{|\alpha|^2}\alpha,
\qquad
x\in\mathbb{R}^N,
\end{gather*}
for the orthogonal projection onto the hyperplane $(\mathbb{R}\alpha)^\perp=\{x,\,  \langle
x,\alpha\rangle=0\}$, so that the ref\/lec\-tion~$s_{\alpha}$ with respect to hyperplane orthogonal to
$\alpha$ is related to $\tau_{\alpha}$ by
\begin{gather*}
\tau_\alpha=\frac{1}{2}(1+s_{\alpha}).
\end{gather*}
The hyperplane $(\mathbb{R}\alpha)^{\perp}$ is the invariant set of $\tau_{\alpha}$. If $\langle
\alpha,\beta\rangle=0$, then the orthogonal projections~$\tau_{\alpha}$ and~$\tau_{\beta}$ commute.
The conjugate of orthogonal projection onto a~hyperplane is again an orthogonal projection onto
a~hyperplane: suppose $u \in O (N)$ and $\alpha \in \mathbb{R}^N\backslash\{0\}$ then
\begin{gather*}
u\tau_{\alpha}u^{-1}=\tau_{u\alpha}.
\end{gather*}
Let $R$ be a~root system and $R'$ a~positive orthogonal subsystem of $R$.
For $\xi \in \mathbb{R}^N$, we def\/ine the dif\/ferential-dif\/ference operator $T_\xi$ by
\begin{gather}
(T_{\xi}f)(x)=\partial_{\xi}f(x)+\sum_{\alpha\in R'}
\kappa(\alpha)\langle\alpha,\xi\rangle\frac{f(x)-f(\tau_{\alpha}x)}{\langle x,\alpha\rangle}.
\label{operat1}
\end{gather}
where $\kappa$ is a~multiplicity function on $R'$.
For $j=1,\dots,N$ denotes~$T_{e_j}$ by~$T_j$.
The operator~$T_\xi$ can be considered as a~deformation of the usual directional derivatives and when
$\kappa = 0$, the operator $T_\xi$ reduces to the corresponding directional derivative.
Furthermore, there is overlap between the notations~\eqref{operat1} and~\eqref{Dunkl}.
In fact, the operator~\eqref{operat1} follows from Dunkl operator after replacing the ref\/lections terms
that exist in~\eqref{Dunkl} by orthogonal projection terms.
\begin{example}\label{example3.1}
In the rank-one case, the root system is of type $A_1$ and the corresponding ref\/lection $s$ and
orthogonal projection $\tau$ are given by
\begin{gather*}
s(x)=-x,
\qquad
\tau(x)=\frac{1}{2}(1+s)(x)=0.
\end{gather*}
The Dunkl-type operator $T_{\kappa}$ associated with the projection $\tau$ and the multiplicity parameters~$\kappa$ $(\kappa \in \mathbb{C})$ is given by
\begin{gather*}
T_{\kappa}f(x)=f'(x)+\kappa \frac{f(x)-f(\tau(x))}{x}=f'(x)+\kappa \frac{f(x)-f(0)}{x}.
\end{gather*}
\end{example}

\begin{example}\label{example3.2}
Let $R=\{\pm(e_1\pm e_2),\pm e_1,\pm e_2\}$ be a~root system of type~$B_2$ in the 2-plane and
$R'=\{e_1\pm e_2\}$ be a~positive orthogonal subsystem in~$R$.
The related Dunkl-type operators to~$R'$ and to the positive parameters $(\kappa_1,\kappa_2)$ are given by
\begin{gather*}
T_1=\partial_x+\kappa_1\frac{f(x,y)-f((x+y)/2,(x+y)/2)}{x-y}+\kappa_2\frac{f(x,y)-f((x-y)/2,(x-y)/2)}{x+y},
\\
T_2=\partial_y-\kappa_1\frac{f(x,y)-f((x+y)/2,(x+y)/2)}{x-y}+\kappa_2\frac{f(x,y)-f((x-y)/2,(x-y)/2)}{x+y}.
\end{gather*}
\end{example}
We denote by $\Pi^N$ the space of polynomials and by $\Pi_n^N$ the subspace of homogenous polynomials of
degree $n.$

Let $R'=\{\alpha_1,\dots,\alpha_n\}$ be a~positive orthogonal subsystem of a~root system $R$.
Consider the operator $\rho_{i}$ def\/ined on $\Pi^N$ by
\begin{gather*}
(\rho_{i}f)(x)=\frac{f(x)-f(\tau_{\alpha_i}x)}{\langle x,\alpha\rangle},
\qquad
i=1,\dots,n.
\end{gather*}
It follows from the equality
\begin{gather*}
(\rho_{j}f)(x)=-\frac{1}{|\alpha_j|^2}\int_0^1\partial_{\alpha_j}
f\left(x-t\frac{\langle x,\alpha_j\rangle}{|\alpha_j|^2}\alpha_j\right)\,{\rm d}t
\end{gather*}
that $T_\xi$ is a~homogeneous operator of degree $-1$ on $\Pi^{N}$, that is, $T_\xi f \in \Pi_{n-1}^{N}$,
for $f \in \Pi_n^{N}$, and leaves $\mathcal{S}(\mathbb{R})$ ($\mathcal{S}(\mathbb{R})$ is the Schwartz space of
rapidly decreasing functions on $\mathbb{R}$) invariant.
\begin{proposition}\label{proposition3.1}
The operators $\rho_i$ $(i=1,\dots,n)$ have the following properties:
\begin{itemize}\itemsep=0pt
\item[$i)$] for $i,j=1,\dots,n$, we have
$
[\rho_i,\rho_j]=0$;

\item[$ii)$] if $\alpha$ is an orthogonal vector to $\alpha_i$, then
$
[\partial_\alpha,\rho_i]=0$,
where the commutator of two operators $A$, $B$ is defined by
$
\left[A,B\right]:=AB-BA$.
\end{itemize}
\end{proposition}

The family $\{\alpha_1,\dots,\alpha_n\}$ is orthogonal, then there exist scalars
$\xi_1,\dots,\xi_n$ and a~vector $\widehat{\xi}\in \mathbb{R}^N$ orthogonal to the subspace
$\mathbb{R}\alpha_1\oplus\dots\oplus\mathbb{R}\alpha_n$ such that
\begin{gather*}
\xi=\sum_{i=1}^n\xi_i\alpha_i+\widehat{\xi}.
\end{gather*}
This allows us to decompose the operator $T_\xi$~\eqref{operat1} associated with $R'$ and the multiplicity
para\-me\-ters $(\kappa_1,\dots,\kappa_n)$ in a~unique way in the form
\begin{gather*}
T_\xi=\sum_{i=1}^n\xi_iT_{\alpha_i}+\partial_{\widehat{\xi}}.
\end{gather*}
We now have all ingredients to state and prove the f\/irst main result of the paper.
\begin{theorem}\label{theorem3.1}
Let $\xi,\eta \in \mathbb{R}^N$, then $[T_\xi,T_\eta]=0.$
\end{theorem}

\begin{proof}
A straightforward computation yields
\begin{gather*}
[T_\xi,T_\eta]=\sum_{i,j=1}^{n}\xi_{i}\eta_{j}[T_{\alpha_i},T_{\alpha_j}]
+[\partial_{\widehat{\xi}},\partial_{\widehat{\eta}}]
+\sum_{i=1}^{n}\xi_{i}[T_{\alpha_i},\partial_{\widehat{\eta}}]-\eta_{i}[T_{\alpha_i},\partial_{\widehat{\xi}}].
\end{gather*}
On the other hand,
\begin{gather*}
[T_{\alpha_i},T_{\alpha_j}]=[\partial_{\alpha_i}+\kappa_i\|\alpha_i\|\rho_i,\partial_{\alpha_j}
+\kappa_j\|\alpha_j\|\rho_j]
\\
\phantom{[T_{\alpha_i},T_{\alpha_j}]}{}
=[\partial_{\alpha_i},\partial_{\alpha_j}]+\kappa_{j}\|\alpha_j\|[\partial_{\alpha_i},\rho_{j}
]-\kappa_{i}\|\alpha_i\|[\partial_{\alpha_j},\rho_{i}]
+\kappa_{i}\kappa_j\|\alpha_i\|\|\alpha_j\|[\rho_{i},\rho_{j}],
\end{gather*}
and
\begin{gather*}
[T_{\alpha_i},\partial_{\widehat{\xi}}]=[\partial_{\alpha_i},\partial_{\widehat{\xi}}
]+\kappa_i\|\alpha_i\|[\rho_i,\partial_{\xi}].
\end{gather*}
From Proposition~\ref{proposition3.1}, we get
\begin{gather*}
[T_{\alpha_i},T_{\alpha_j}]=0
\qquad
\text{and}
\qquad
[T_{\alpha_i},\partial_{\widehat{\xi}}]=0.
\end{gather*}
This proves the result.
\end{proof}

One important consequence of the Theorem~\ref{theorem3.1}, is that the operators $T_{\alpha_1},\ldots,T_{\alpha_m}$
generate a~commutative algebra.

\section{Intertwining operator}\label{Sect4}

In this section, we give an intertwining operator between $T_\xi$ and the directional derivative $\partial_\xi$.
Consider a~positive orthogonal subsystem $R'=\{\alpha_1,\ldots,\alpha_n\}$
composed of~$n$ vectors in a root system $R$,
and $\kappa=(\kappa_1,\ldots,\kappa_n)\in \mathbb{C}^n $ and $\xi \in \mathbb{R}^N$.
The associated Dunkl-type operator~$T_\xi$ with~$R'$ and~$\kappa$ takes the form
\begin{gather*}
(T_{\xi}f)(x)=\partial_{\xi}f(x)+\sum_{j=1}
^n\kappa_j\langle\alpha_j,\xi\rangle\frac{f(x)-f(\tau_{\alpha_j}x)}{\langle x,\alpha_j\rangle}.
\end{gather*}
Let $h: \mathbb{R}^n\times\mathbb{R}^N\rightarrow \mathbb{R}^N$ be the function def\/ined by
\begin{gather*}
h(t,x)=x+\sum_{j=1}^n(t_j-1)\frac{\langle x,\alpha_j\rangle}{|\alpha_j|^2}\alpha_j,
\end{gather*}
where $t=(t_1,\ldots,t_n)\in \mathbb{R}^n$ and $x\in \mathbb{R}^N.$

We def\/ine
\begin{gather}
\chi_\kappa(f)(x)=\frac{1}{\Gamma(\kappa)}\int_{[0,1]^n}f(h(t,x))w(t)\,{\rm d}t ,
\label{inter}
\end{gather}
where $w(t)=\prod\limits_{j=1}^n(1-t_j)^{\kappa_j-1}$ and
$\Gamma(\kappa)=\prod\limits_{j=1}^n\Gamma(\kappa_j).$
\begin{theorem}\label{theorem4.1}
Let $f\in C^{\infty}(\mathbb{R}^N)$, then we have
\begin{gather*}
T_\xi\circ\chi_\kappa f(x)=\chi_\kappa\circ\partial_\xi f(x).
\end{gather*}
\end{theorem}
\begin{proof}
For $j=1,\ldots,n$, we denote by $\theta_j $ the orthogonal projection in $\mathbb{R}^n$ with respect to
the hyperplane $(\mathbb{R}e_j)^\perp$ orthogonal to the vector $e_j$ of the canonical basis $(e_1,\ldots,
e_n)$ of $\mathbb{R}^n$. The orthogonal projection $\theta_j$ acts on $\mathbb{R}^n$ as
\begin{gather*}
\theta_j(t)=(t_1,\ldots,t_{j-1},0,t_{j+1},\ldots,t_n).
\end{gather*}
The system $R$ is orthogonal, then for $j=1,\ldots, n$, we have
\begin{gather*}
h(t,\tau_{\alpha_j}x)=\tau_{\alpha_j}x+\sum_{k=1}^n(t_k-1)\frac{\langle\tau_{\alpha_j}x,\alpha_k\rangle}
{|\alpha_k|^2}\alpha_k
\\
\phantom{h(t,\tau_{\alpha_j}x)}{}
=x-\frac{\langle x,\alpha_j\rangle}{|\alpha_j|^2}\alpha_j+\sum_{k=1,k\neq j}
^n(t_k-1)\frac{\langle x,\alpha_k\rangle}{|\alpha_k|^2}\alpha_k
=h(\theta_jt,x).
\end{gather*}
Let $f \in C^\infty(\mathbb{R}^N)$ and $\xi \in \mathbb{R}^N$. The mapping $x\rightarrow h(t,x)$ is linear
on $\mathbb{R}^N$, then we can write
\begin{gather*}
\partial_\xi(f(h(t,x)))=\partial_{h(t,\xi)}f(h(t,x))
=\partial_{\xi}f(h(t,x))+\sum_{j=1}^n(t_j-1)\frac{\langle\xi,\alpha_j\rangle}{|\alpha_j|^2}
\partial_{\alpha_j}f(h(t,x)).
\end{gather*}
Hence,
\begin{gather*}
\partial_\xi\chi_\kappa(f)(x)=\frac{1}{\Gamma(\kappa)}\int_{[0,1]^n}\partial_{\xi}(f(h(t,x)))w(t)\,{\rm d}t
=\frac{1}{\Gamma(\kappa)}\int_{[0,1]^n}\partial_{\xi}f(h(t,x))w(t)\,{\rm d}t
\\
\phantom{\partial_\xi\chi_\kappa(f)(x)=}{}
+\frac{1}{\Gamma(\kappa)}\sum_{j=1}^n\frac{\langle\xi,\alpha_j\rangle}{|\alpha_j|^2}\int_{[0,1]^n}
(t_j-1)\partial_{\alpha_j}f(h(t,x))w(t)\,{\rm d}t .
\end{gather*}
Since we can write
\begin{gather*}
\partial_{t_j}f(h(t,x))=\frac{\langle x,\alpha_k\rangle}{|\alpha_k|^2}\partial_{\alpha_j}f(h(t,x))
\end{gather*}
and
\begin{gather*}
\int_0^1(1-t_j)^{\kappa_j}\partial_{t_j}
f(h(t,x))\,{\rm d}t =-f(h(\theta_j(t),x))+\kappa_j\int_0^1(1-t_j)^{\kappa_j-1}f(h(t,x))\,{\rm d}t,
\end{gather*}
we are lead to
\begin{gather*}
\int_{[0,1]^n}\partial_{\alpha_j}f(h(t,x))(t_j-1)w(t)\,{\rm d}t =\frac{|\alpha_j|^2}{\langle x,\alpha_j\rangle}
\int_{[0,1]^n}\partial_{t_j}f(h(t,x))(t_j-1)w(t)\,{\rm d}t
\\
\qquad{}
=\kappa_j\frac{|\alpha_j|^2}{\langle x,\alpha_j\rangle}\int_{[0,1]^n}(f(h(\theta_j(t),x))-f(h(t,x)))w(t)\,{\rm d}t
\\
\qquad{}
=-\kappa_j\Gamma(\kappa)\frac{|\alpha_j|^2}{\langle x,\alpha_j\rangle}
\big(\chi_\kappa(f)(x)-\chi_\kappa(f)(\tau_{\alpha_j}x)\big).
\end{gather*}
This, combined with the last expression of $\partial_\xi(\chi_\kappa f)(x)$, yields
\begin{gather*}
\partial_\xi\chi_\kappa(f)(x)=\chi_{\kappa}(\partial_{\xi}f)(x)-\sum_{j=1}
^n\kappa_j\langle\xi,\alpha_j\rangle\frac{\chi_\kappa(f)(x)-\chi_\kappa(f)(\tau_jx)}
{\langle x,\alpha_j\rangle}.
\end{gather*}
Therefore,
\begin{gather*}
T_\xi(\chi_\kappa f)(x)=\chi_{\kappa}(\partial_\xi f)(x).\tag*{\qed}
\end{gather*}
\renewcommand{\qed}{}
\end{proof}

\section{The one variable case}\label{Sect5}

The specialization of this theory to the one variable case has its own
interest, because everything can be done there in a~much more explicit way and new results for special
functions in one variable can be obtained.
In this setting there is only one Dunkl-type operator $T_{\kappa}$ associated up to scaling and it equals to
\begin{gather}
T_{\kappa}f(x)=f'(x)+\kappa \frac{f(x)-f(0)}{x}.
\label{Di}
\end{gather}
This operator leaves the space of polynomials invariant and acts on the monomials as
\begin{gather*}
T_{\kappa}1=0,
\qquad
T_{\kappa}x^n=(n+\kappa)x^{n-1},
\qquad
n=1,2,\dots.
\end{gather*}
Its square is given by
\begin{gather*}
T^2_{\kappa}f(x)=f''(x)+ \frac{2\kappa}{x}
f'(x)+\kappa(\kappa-1) \frac{f(x)-f(0)}{x^2}- \frac{\kappa(\kappa+1)}{x}f'(0).
\end{gather*}
Consider the conf\/luent hypergeometric function (see~\cite[\S~7.1]{Temme})
\begin{gather*}
M(a,b;z)=\sum_{n=0}^{\infty}\frac{(a)_n}{(b)_n}\frac{z^n}{n!},
\end{gather*}
where $(a)_n$ is the Pochhammer symbol def\/ined by
\begin{gather*}
(a)_n=\frac{\Gamma(a+n)}{\Gamma(a)}.
\end{gather*}
This is a~solution of the conf\/luent hypergeometric dif\/ferential equation
\begin{gather*}
zy''(z)+(b-z)y'(z)=ay(z).
\end{gather*}
This function possesses the following Poisson integral representation (see~\cite[\S~7.1]{Temme})
\begin{gather}
M(a,b;z)=\frac{\Gamma(b)}{\Gamma(a)\Gamma(b-a)}\int_{0}^{1}t^{a-1}(1-t)^{b-a-1}e^{zt}\,{\rm d}t,
\qquad
\Re(b)>\Re(a)>0.\label{int1}
\end{gather}
\begin{theorem}\label{theorem5.1}
For $\lambda \in \mathbb{C}$ and $\kappa>-1$, the problem
\begin{gather}
T_{\kappa}f(x)=i\lambda f(x),
\qquad
f(0)=1,
\label{eigen}
\end{gather}
has a~unique analytic solution $M_{\kappa}(i\lambda x)$ given by
\begin{gather}
M_{\kappa}(i\lambda x)=M(1,\kappa+1;i\lambda x).
\label{kummer2}
\end{gather}
\end{theorem}
\begin{proof}
Searching a~solution of~\eqref{eigen} in the form $f(z)=\sum\limits_{n=0}^{\infty}a_nx^n$. Replacing
in~\eqref{eigen}, we obtain
\begin{gather*}
\sum_{n=0}^{\infty}(n+1+\kappa)a_{n+1}x^n=i\lambda\sum_{n=0}^{\infty}a_nx^n.
\end{gather*}
Thus,
\begin{gather*}
a_{n+1}=\frac{i\lambda}{n+1+\kappa}a_n
\qquad \text{and}\qquad
a_n=\frac{(i\lambda)^n}{(\kappa+1)_n}.\tag*{\qed}
\end{gather*}
\renewcommand{\qed}{}
\end{proof}

\begin{remark}\label{remark5.1}
Multiply the equation~\eqref{eigen} by $x$ and dif\/ferentiating both sides, we see that a~function~$u$ of
class $C^2$ on $\mathbb{R}$, is a~solution of the equation~\eqref{eigen}, if and only if, it is a~solution
of the generalized eigenvalue problem
\begin{gather*}
xu''+(\kappa+1)u'=i\lambda(xu'+u).
\end{gather*}
\end{remark}

\begin{proposition}\label{proposition5.1}
The function $\mathbf{M}_\kappa(z)$ defined by
\begin{gather}
\mathbf{M}_\kappa(z)=\frac{M_\kappa(z)}{\Gamma(\kappa+1)}=\sum_{n=0}^{\infty}\frac{z^n}{\Gamma(\kappa+1+n)}
\label{Kummer1}
\end{gather}
satisfies the following properties:
\begin{itemize}\itemsep=0pt
\item[$(i)$] $\mathbf{M}_\kappa(z)$ is analytic in $\kappa$ and $z$;
\item[$(ii)$] $\mathbf{M}_0(z)=e^z$;
 \item[$(iii)$] for $\Re(\kappa)>0$, the function $\mathbf{M}_\kappa(z)$,
possesses the integral representation
\begin{gather*}
\mathbf{M}_\kappa(z)= \frac{1}{\Gamma(\kappa)}\int_{0}^{1}(1-t)^{\kappa-1}e^{zt}\,{\rm d}t;
\end{gather*}
\item[$(iv)$] for $\Re(\kappa)>0$, we have
\begin{gather*}
\big|\mathbf{M}^{(n)}_\kappa(z)\big|\leq|z|^ne^{\Re(z)},
\qquad
n\in\mathbb{N}, \qquad z\in\mathbb{C},
\end{gather*}
in particular,
\begin{gather*}
|\mathbf{M}_\kappa(i\lambda x)|\leq1,
\qquad
\lambda,x\in\mathbb{R};
\end{gather*}
\item[$(v)$] for $\Re(\kappa)>0$, and all $x\in \mathbb{R}^\ast$,
\begin{gather*}
 \lim_{\lambda\rightarrow+\infty}\mathbf{M}_\kappa(i\lambda x)=0.
\end{gather*}
\end{itemize}
\end{proposition}

\begin{proof}
 (i) and (ii) are immediate.
(iii) follows from~\eqref{int1}.
For $n\in \mathbb{N}$, we have
\begin{gather*}
\mathbf{M}^{(n)}_\kappa(z)=\frac{z^n}{\Gamma(\kappa)}\int_0^1(1-t)^{\kappa}t^ne^{zt}\,{\rm d}t.
\end{gather*}
So we f\/ind
\begin{gather*}
|\mathbf{M}^{(n)}_{\kappa}(z)|\leq\frac{|z|^n}{\Gamma(\kappa)}\int_0^1(1-t)^{\kappa}e^{\Re(z)t}\,{\rm d}t
\leq|z|^ne^{\Re(z)}.
\end{gather*}
This proves (iv).
(v) follows from (iii) and the Riemann--Lebesgue lemma.
\end{proof}

\begin{definition} \label{definition5.1}
We def\/ine the Kummer transform on $L^1(\mathbb{R})$ by
\begin{gather*}
\forall\, \lambda\in\mathbb{R},
\qquad
\mathcal{F}_{\kappa}(f)(\lambda)=\int_{\mathbb{R}}f(x)\mathbf{M}_\kappa(i\lambda x)(x)\,{\rm d}x .
\end{gather*}
When $\kappa=0$, the transformation $\mathcal{F}_{0}$ reduces to the usual Fourier transform $\mathcal{F}$
that is given by
\begin{gather*}
\mathcal{F}(f)(\lambda)=\int_{\mathbb{R}}f(x)e^{i\lambda x}\,{\rm d}x .
\end{gather*}
\end{definition}

\begin{theorem}  \label{theorem5.2}
Let $f$ be a~function in $L^1(\mathbb{R})$ then $\mathcal{F}_{\kappa}(f)$ belongs to $C_0(\mathbb{R})$,
where $C_0(\mathbb{R})$ is the space of continuous functions having zero as limit at the infinity.
Furthermore,
\begin{gather*}
\|\mathcal{F}_{\kappa}(f)\|_\infty\leq\|f\|_{1}.
\end{gather*}
\end{theorem}

\begin{proof}
It's clear that $\mathcal{F}_{\kappa}(f)$ is a~continuous function on $\mathbb{R}$.
From Proposition~\ref{proposition5.1}, we get for all $x\in \mathbb{R}^\ast$,
\begin{gather*}
 \lim_{\lambda\rightarrow\infty}f(x)\mathbf{M}_{\kappa}(i\lambda x)=0
\qquad \text{and}\qquad
 |f(x)\mathbf{M}_{\kappa}(i\lambda x)|\leq|f(x)|.
\end{gather*}
Since $f$ is in $L^1(\mathbb{R})$, we conclude by using the dominated convergence theorem that
$\mathcal{F}_{\kappa}(f)$ belongs to $C_0(\mathbb{R})$ and
\begin{gather*}
\|\mathcal{F}_{\kappa}(f)\|_\infty\leq\|f\|_{1}.
\end{gather*}
We now turn to exhibit a~relationship between the Kummer transform and the Fourier transform.
The crucial idea is to use the intertwining operator $\chi_\kappa$.
We denote by $C^{\infty}(\mathbb{R})$ the space of inf\/initely dif\/ferentiable functions $f$ on
$\mathbb{R}$, provided with the topology def\/ined by the semi norms
\begin{gather*}
\|f\|_{n,a}=\sup_{\substack{0\leq k\leq n\\ x\in[-a,a]}}\big|f^{(k)}(x)\big|,
\qquad
a>0,\qquad
n\in\mathbb{N}.
\end{gather*}
In the rank-one case the intertwining operator~\eqref{inter} becomes
\begin{gather}
(\chi_\kappa f)(x)=\frac{1}{\Gamma(\kappa)}\int_0^1(1-t)^{\kappa-1}f(tx)\,{\rm d}t .
\label{intert4}
\end{gather}
This operator is a~particular case of the so called Erd\'{e}lyi--Kober fractional integral
$I^{\gamma,\delta}$, which is given by (see~\cite{Kob})
\begin{gather*}
(I^{\gamma,\delta}f)(x)=\frac{1}{\Gamma(\delta)}\int_0^1(1-t)^{\delta-1}t^{\gamma}f(tx)\,{\rm d}t ,
\qquad
\delta>0,
\qquad
\gamma\in\mathbb{R}.
\end{gather*}
It was shown in~\cite[\S~3]{Luchko}, that the Erd\'{e}lyi--Kober fractional integral has a~left-inverse
\begin{gather}
D^{\gamma,\delta}I^{\gamma,\delta}f=f,
\qquad
f\in C^\infty(\mathbb{R}),\label{kob1}
\end{gather}
where
\begin{gather*}
D^{\gamma,\delta}=\prod_{k=1}^n\left(\gamma+k+x\frac{{\rm d}}{{\rm d}x }\right)I^{\gamma+\delta,n-\delta},
\end{gather*}
and $n=\lceil\delta\rceil$ ($\lceil\delta\rceil$ denotes the ceiling function the smallest integer $\geq\delta$).

As a~consequence of Theorem~\ref{theorem4.1}, we deduce that the operator $\chi_\kappa$~\eqref{intert4} has the
fundamental intertwining property
\begin{gather*}
T_\kappa\circ\chi_\kappa=\chi_\kappa\circ\frac{{\rm d}}{{\rm d}x }.
\end{gather*}
We regard it as a~second main result since it allows us to move from the complicated operator~$T_\kappa$
def\/ined in~\eqref{Di} to the simple derivative operator $\frac{{\rm d}}{{\rm d}x }$.
\end{proof}

\begin{theorem}\label{theorem5.3}
Let $\kappa>0$, the operator $\chi_\kappa$ is a~topological isomorphism from $C^{\infty}( \mathbb{R})$ onto
itself and its inverse $\chi_\kappa^{-1}$ is given for all $f\in C^{\infty}( \mathbb{R})$ by
\begin{gather*}
\chi_\kappa^{-1}f(x)=D^{0,\kappa}f(x)=\prod_{j=1}^n\left(j+x\frac{{\rm d}}{{\rm d}x }\right)(I^{\kappa+1,n-\kappa}f)(x),
\end{gather*}
where $n=\lceil\kappa\rceil.$
\end{theorem}

\begin{proof}
Let $a>0$ and $f\in C^\infty( \mathbb{R})$.
For $x\in[0,a]$, $t\in [0,1]$ and $l\in \mathbb{N}$, we have the following estimate
\begin{gather*}
|t^l(1-t)^{\kappa-1}f^{(l)}(xt)|\leq\|f\|_{l,a}(1-t)^{\kappa-1}
\qquad
\text{and}
\qquad
\int_0^1(1-t)^{\kappa-1}\,{\rm d}t =\frac{1}{\kappa}.
\end{gather*}
By the theorem of derivation under the integral sign, we can prove that
\begin{gather*}
\chi_{\kappa}f\in C^\infty(\mathbb{R})
\qquad
\text{and}
\qquad
\|\chi_\kappa(f)\|_{l,a}\leq\frac{1}{\Gamma(\kappa+1)}\|f\|_{l,a}.
\end{gather*}
Then $\chi_\kappa$ is a~linear continuous mapping from $C^\infty(\mathbb{R})$ onto its self.
From formula~\eqref{kob1} the operator
\begin{gather*}
D^{0,\kappa}=\prod_{j=1}^n\left(j+x\frac{{\rm d}}{{\rm d}x }\right)\circ I^{\kappa+1,n-\kappa}
\end{gather*}
is a~left-inverse of $\chi_{\kappa}$.
This shows that $\chi_\kappa$ is injective and $D^{0,\kappa}$ is surjective.
So it suf\/f\/ices to prove that $D^{0,\kappa}$ is injective.

Let $f$ be a~function in $C^\infty(\mathbb{R})$ such that
$ 
D^{0,\kappa}f=0$.
Then the function $g=I^{\kappa+1,n-\kappa}f\in C^{\infty}(\mathbb{R})$ is a~solution of the linear
dif\/ferential equation
\begin{gather*}
\prod_{j=1}^n\left(1+j+x\frac{{\rm d}}{{\rm d}x }\right)y(x)=0.
\end{gather*}
Since, the last dif\/ferential equation has a~unique $C^\infty$-solution, which is equal to $y(x)=0$, it
follows that $g=0$.

From~\eqref{kob1} the operator $I^{\kappa+1,\kappa}$ has a~left-inverse, then $f=0$.
This shows that $\chi_{\kappa}$ is a~bijective operator.
\end{proof}

Let $\kappa>0$, we def\/ine the dual intertwining operator ${}^t\chi_\kappa$ on $\mathcal{D}(\mathbb{R})$
($\mathcal{D}(\mathbb{R})$ is the space of $C^\infty$-functions on $\mathbb{R}$ with compact support) by
\begin{gather*}
\big({}^t\chi_\kappa f\big)(x)=\frac{1}{\Gamma(\kappa)}\int_{|x|}^{+\infty}(t-|x|)^{\kappa-1}t^{-\kappa}
f({\rm sgn}(x)t)\,{\rm d}t ,
\qquad
x\in\mathbb{R}\setminus\{0\}.
\end{gather*}
\begin{proposition}\label{proposition5.2}
The operator ${}^t\chi_\kappa$ is a~topological automorphism of $\mathcal{D}(\mathbb{R})$, and satisfies
the transmutation relation:
\begin{gather*}
\int_{\mathbb{R}}(\chi_{\kappa}f)(x)g(x)\,{\rm d}x =\int_{\mathbb{R}}f(x)\big({}^t\chi_{\kappa}g\big)(x)\,{\rm d}x ,
\qquad
f\in C^\infty(\mathbb{R}).
\end{gather*}
\end{proposition}

\begin{proof}
 Let $f\in C^{\infty}(\mathbb{R})$ and $g \in \mathcal{D}(\mathbb{R})$, we have
\begin{gather*}
\int_{\mathbb{R}}(\chi_{\kappa}f)(x)g(x)\,{\rm d}x =\frac{1}{\Gamma(\kappa)}\int_0^{+\infty}
\int_0^x(x-t)^{\kappa-1}f(t)\,{\rm d}t g(x)x^{-\kappa}\,{\rm d}x
\\
\phantom{\int_{\mathbb{R}}(\chi_{\kappa}f)(x)g(x)\,{\rm d}x =}{}
-\frac{1}{\Gamma(\kappa)}\int_0^{\infty}\int_0^x(x-t)^{\kappa-1}f(-t)\,{\rm d}t g(-x)x^{-\kappa}\,{\rm d}x .
\end{gather*}
Using Fubini's theorem and a~change of variable, we get
\begin{gather*}
\int_{\mathbb{R}}(\chi_{\kappa}f)(x)g(x)\,{\rm d}x =\frac{1}{\Gamma(\kappa)}\int_0^{+\infty}
\int_t^\infty x^{-\kappa}(x-t)^{\kappa-1}g(x)\,{\rm d}x f(t)\,{\rm d}t
\\
\phantom{\int_{\mathbb{R}}(\chi_{\kappa}f)(x)g(x)\,{\rm d}x =}{}
+\frac{1}{\Gamma(\kappa)}\int_{-\infty}^0\int_{-t}^\infty x^{-\kappa}(x+t)^{\kappa-1}g(-x)\,{\rm d}x f(t)\,{\rm d}t .
\end{gather*}
Therefore,
\begin{gather*}
\int_{\mathbb{R}}(\chi_{\kappa}f)(x)g(x)\,{\rm d}x =\frac{1}{\Gamma(\kappa)}\int_{\mathbb{R}}\int_{|t|}
^\infty x^{-\kappa}(x-|t|)^{\kappa-1}g(sign(t)x)\,{\rm d}x f(t)\,{\rm d}t
\\
\phantom{\int_{\mathbb{R}}(\chi_{\kappa}f)(x)g(x)\,{\rm d}x }{}
=\int_{\mathbb{R}}f(t)\big({}^t\chi_{\kappa}g\big)(t)\,{\rm d}t .\tag*{\qed}
\end{gather*}
\renewcommand{\qed}{}
\end{proof}

\begin{proposition}\label{proposition5.3}
Let $\kappa>0$, the Kummer transform $\mathcal{F}_\kappa$ satisfies the decomposition
\begin{gather*}
\mathcal{F}_{\kappa}(f)=\mathcal{F}\circ{}^t\chi_\kappa(f),
\qquad
f\in\mathcal{D}(\mathbb{R}).
\end{gather*}
\end{proposition}
\begin{proof}
The result follows from Proposition~\ref{proposition5.2}.
\end{proof}

\section{Multivariable case}\label{Sect6}

\subsection{Direct product setting}

In this subsection, we consider the direct product of the  one-dimensional models, which means that the Weyl group
of the corresponding subsystem of root system is a subgroup of $\mathbb Z_2^N$.

We denote by $\tau_k$ (for each $k$ from $1$ to $N$)
the orthogonal projection with respect to the hyperplane orthogonal to $e_k$,
that is to say for every $x=(x_1,\ldots,x_N) \in \mathbb R^N$
\begin{gather*}
\tau_k(x)=x-\frac{\left\langle x,e_k\right\rangle}{|e_k|^2}e_k=(x_1,\ldots,x_{k-1},0,x_{k+1},\ldots,x_N).
\end{gather*}
Let $\kappa=(\kappa_1,\kappa_2,\ldots,\kappa_N) \in \mathbb{C}^N$. The associated Dunkl type operators
$T_j$ for $j=1,\dots,N$, are given for $x \in \mathbb R^N$ by
\begin{gather*}
T_j f(x)=\partial_jf(x)+\sum_{l=1}^N\kappa_l\frac{f(x)-f(\tau_l(x))}
{\left\langle x,e_l\right\rangle}\left\langle e_k,e_l\right\rangle
\\
\phantom{T_j f(x)}{}
=\partial_jf(x)+\kappa_j\frac{f(x)-f(x_1,\ldots,x_{j-1},0,x_{j+1},\ldots,x_N)}{x_j}.
\end{gather*}
These operators form a~commuting system.
The generalized Laplacian associated with $T_j$ is def\/ined in a~natural way as
\begin{gather*}
\Delta_\kappa=\sum_{j=1}^NT^2_j.
\end{gather*}
A straightforward computation yields
\begin{gather*}
\Delta_{\kappa}=\Delta+2\sum_{j=1}^N\kappa_j x^{-1}_j\partial_jf(x)-\sum_{j=1}^N(\kappa^2_j+\kappa_j)x^{-1}
_j\partial_jf(x_1,\ldots,x_{j-1},0,x_{j+1},\ldots,x_N)
\\
\phantom{\Delta_{\kappa}=}{}
+\sum_{j=1}^N(\kappa^2_j-\kappa_j)x^{-2}_j\big(f(x)-f(x_1,\ldots,x_{j-1},0,x_{j+1},\ldots,x_N)\big).
\end{gather*}

This operator will play in our context a~similar role to that of the Euclidean Laplacian in the classical
harmonic analysis.
Obviously, the trivial choice of the multiplicity function $\kappa=0$, reduces our situation to the
analysis related to the classical Laplacian $\Delta$.

Let $\kappa=(\kappa_1,\ldots,\kappa_N) \in (0,\infty)^N$. For $x,\lambda \in \mathbb{R}^N$, we consider
the function $M_\kappa(\lambda,x)$ which is given as the tensor products
\begin{gather*}
M_\kappa(\lambda,x)=\prod_{j=1}^NM_{\kappa_j}(i\lambda_jx_j).
\end{gather*}
\begin{theorem}\label{theorem6.1}
For $\lambda=(\lambda_1,\ldots,\lambda_N) \in \mathbb{C}^N$, the function $M_\kappa(\lambda,x)$ is the
unique analytic solution of the system
\begin{gather*}
T_{\xi}u(x)=i\langle\lambda,\xi\rangle u(x),
\qquad
u(0)=1,
\qquad
\forall\, \xi\in\mathbb{C}^N.
\end{gather*}
\end{theorem}

\subsection[Dunkl-type operators associated to an orthogonal subsystem in a~root system of type $A_{N-1}$]
{Dunkl-type operators associated to an orthogonal subsystem\\ in a~root system of type $\boldsymbol{A_{N-1}}$}

Let $R$ be a~root system of type $A_{N-1}$
\begin{gather*}
R=\{\pm(e_i-e_j),\, 1\leq i<j\leq N\}.
\end{gather*}
Def\/ine a~positive orthogonal subsystem $R'=\{\alpha_1,\dots,\alpha_{[N/2]}\}$ of $R$ by setting:
\begin{gather*}
\alpha_i=e_{2i-1}-e_{2i},
\qquad
i=1,\dots,[N/2].
\end{gather*}
We denote by $\tau_j$ (for each $j$ from $1$ to $[N/2]$) the orthogonal projection onto the hyperplane
perpendicular to $\alpha_j$, that is to say for every $x = (x_1, \dots, x_N)\in \mathbb{R}^N$
\begin{gather*}
\tau_ix=(x_1,\dots,\overline{x}_{2i-1},\overline{x}_{2i},\dots,x_N),
\end{gather*}
where
$
\overline{x}_{2i-1}=\overline{x}_{2i}=\frac{1}{2}(x_{2i-1}+x_{2i})$,
$i=1,\dots,[N/2]$.
The vector $\xi\in\mathbb{R}^N$ can be decomposed uniquely in the form{\samepage
\begin{gather*}
\xi=\sum_{i=1}^{[N/2]}\xi_i(e_{2i-1}-e_{2i})+\widehat{\xi},
\end{gather*}
where $\widehat{\xi}$ is an orthogonal vector to the linear space generated by
$R'=\{\alpha_1,\dots,\alpha_{[N/2]}\}$.}

A straightforward computation shows that the operator $T_\xi$ $(\xi\in\mathbb{R}^N)$ associated with $R'$
and the multiplicity parameters $(\kappa_1,\dots, \kappa_{[N/2]})$ has the following decomposition
\begin{gather*}
T_\xi=\sum_{i=1}^{[N/2]}\xi_{i}T_{\alpha_i}+\partial_{\widehat{\xi}}=\sum_{i=1}^{2[N/2]}(-1)^{i+1}
\xi_{[\frac{i+1}{2}]}T_{i}+\partial_{\widehat{\xi}},
\end{gather*}
where
\begin{gather*}
T_{i}=\partial_{i}-(-1)^{i}\kappa_{[\frac{i+1}{2}]}\rho_{[\frac{i+1}{2}]},
\qquad
i=1,\dots,2[N/2],
\end{gather*}
and
\begin{gather*}
(\rho_{i}f)(x)=\frac{f(x)-f(\tau_ix)}{x_{2i-1}-x_{2i}}.
\end{gather*}
The intertwining operator~\eqref{inter} becomes
\begin{gather*}
\chi_\kappa(f)(x)=\frac{1}{\Gamma(\kappa)}\int_{[0,1]^n}f(h(t,x))w(t)\,{\rm d}t ,
\end{gather*}
where
\begin{gather*}
h(t,x)=x+\sum_{i=1}^{[N/2]}\frac{t_i-1}{2}(x_{2i-1}-x_{2i})(e_{2i-1}-e_{2i}).
\end{gather*}
\begin{proposition}\label{proposition6.1}
Let $\lambda=(\lambda_1, \dots,\lambda_N)\in\mathbb{C}^N$, and
$\kappa=(\kappa_1,\dots,\kappa_{[N/2]})\in (0,\infty)^{[N/2]}$. The following eigenvalue problem
\begin{gather}
T_\xi f=i\langle\lambda,\xi\rangle f,
\qquad
f(0)=1,
\qquad
\forall\,\xi\in\mathbb{C}^N,
\label{eigen33}
\end{gather}
has a~unique analytic solution $M_{\kappa}(\lambda,x)$ given by
\begin{gather*}
M_{\kappa}(\lambda,x)=e^{i\langle\lambda,h(0,x)\rangle}\prod_{j=1}^{[N/2]}M_{\kappa_j}
\left(\frac{i}{2}(\lambda_{2j-1}-\lambda_{2j})(x_{2j-1}-x_{2j})\right).
\end{gather*}
\end{proposition}
\begin{proof}
According to Theorem~\ref{theorem4.1}, $\chi_\kappa$ is an intertwining operator between~$T_\xi$ and~$\partial_\xi$. So,
the function $\chi_\kappa(e^{i\langle\lambda,\cdot\rangle})$ is the unique $C^{\infty}$-solution of
problem~\eqref{eigen33}.

Since we can write
\begin{gather*}
\langle\lambda,h(t,x)\rangle=\langle\lambda,h(0,x)\rangle+\sum_{j=1}^{[N/2]}\frac{t_j}{2}(\lambda_{2j-1}
-\lambda_{2j})(x_{2j-1}-x_{2j}),
\end{gather*}
we are lead to
\begin{gather*}
M_{\kappa}(\lambda,x)=\frac{e^{i\langle\lambda,h(0,x)\rangle}}{\Gamma(\kappa)}\int_{[0,1]^{[N/2]}}
e^{\frac{i}{2}\sum\limits_{j=1}^{[N/2]}t_j(\lambda_{2j-1}-\lambda_{2j})(x_{2j-1}-x_{2j})}w(t)\,{\rm d}t
\\
\phantom{M_{\kappa}(\lambda,x)}{}
=e^{i\langle\lambda,h(0,x)\rangle}\prod_{j=1}^{[N/2]}\frac{1}{\Gamma(\kappa_j)}\int_{[0,1]}e^{\frac{i}{2}
(\lambda_{2j-1}-\lambda_{2j})(x_{2j-1}-x_{2j})}(1-t_j)^{\kappa_j-1}\,{\rm d}t_j.
\end{gather*}
If we now use~\eqref{kummer2} and~\eqref{Kummer1} we get
\begin{gather*}
M_{\kappa}(\lambda,x)=e^{i\langle\lambda,h(0,x)\rangle}\prod_{j=1}^{[N/2]}M_{\kappa_j}
\left(\frac{i}{2}(\lambda_{2j-1}-\lambda_{2j})(x_{2j-1}-x_{2j})\right).\tag*{\qed}
\end{gather*}
\renewcommand{\qed}{}
\end{proof}

\subsection[Dunkl-type operators associated to orthogonal subsystem in root system of type
$B_N$]{Dunkl-type operators associated to orthogonal subsystem\\
 in root system of type $\boldsymbol{B_N}$}

Throughout this subsection $R$ is a~root system of type $B_{N}$ which is given by
\begin{gather*}
R=\{\pm e_i\pm e_j,
\,
 1\leq i<j\leq N;\;
 \pm e_i
\,
 1\leq i\leq N\},
\end{gather*}
and $R'$ is a~positive orthogonal subsystem $R'$ in the root system $R$ given by
\begin{gather*}
R'=\{\alpha^\pm_i=e_{2i-1}\pm e_{2i},
\,
1\leq i\leq[N/2]\}.
\end{gather*}
Denote by $\tau^{\pm}_i$ (for each $i$ from $1$ to $[N/2]$) the orthogonal projection onto the hyperplane
perpendicular to $\alpha^{\pm}_i$, that is to say for every $x = (x_1, \dots, x_N)\in \mathbb{R}^N$
\begin{gather*}
\tau^{\pm}_ix=\big(x_1,\dots,\overline{x}^{\pm}_{2i-1},\overline{x}^{\pm}_{2i},\dots,x_N\big),
\end{gather*}
where
$
\overline{x}^{\pm}_{2i-1}=\overline{x}^{\pm}_{2i}=\frac{1}{2}(x_{2i-1}\pm x_{2i})$.
In this case, the Dunkl type operator $T_{\xi}$ associated with $R'$ and the multiplicity parameters
$\big(\kappa_1^\pm,\dots,\kappa^\pm_{[N/2]}\big)$ takes the form
\begin{gather*}
(T_{\xi}f)(x)=\partial_{\xi}f(x)+\sum_{j=1}^{[N/2]}
\kappa^-_j\langle\alpha^-_j,\xi\rangle\frac{f(x)-f(\tau^-_jx)}{\langle x,\alpha^-_j\rangle}
+\kappa^+_j\langle\alpha^+_j,\xi\rangle\frac{f(x)-f(\tau^+_jx)}{\langle x,\alpha^+_j\rangle}.
\end{gather*}
In particular, for $i=1, \dots, 2[N/2]$ we have
\begin{gather*}
T_{i}=\partial_{i}-(-1)^{i}\kappa^{-}_{[\frac{i+1}{2}]}\rho^{-}_{[\frac{i+1}{2}]}
+\kappa^{+}_{[\frac{i+1}{2}]}\rho^{+}_{[\frac{i+1}{2}]}.
\end{gather*}
where
\begin{gather*}
(\rho^{\pm}_{i}f)(x)=\frac{f(x)-f(\tau^{\pm}_ix)}{x_{2i-1}\pm x_{2i}}.
\end{gather*}
The operator $T_\xi$ has also the following decomposition
\begin{gather*}
T_\xi=\sum_{i=1}^{2[N/2]}\Big(\xi^+_{[\frac{i+1}{2}]}+(-1)^{i+1}\xi^-_{[\frac{i+1}{2}]}\Big)T_{i}
+\varepsilon\xi_N\partial_{N},
\end{gather*}
where
\begin{gather*}
\xi=\sum_{i=1}^{[N/2]}\xi_{i}^{+}\alpha_i^{+}+\xi_{i}^{-}\alpha_i^{-}+\varepsilon\xi_N e_N,
\qquad
\text{and}
\qquad
\varepsilon=
\begin{cases}
1
&
\text{if $N$ is odd},
\\
0
&
\text{if $N$ is even}.
\end{cases}
\end{gather*}
\begin{proposition}\label{proposition6.2}
Let $\lambda=(\lambda_1,\dots, \lambda_N)\in\mathbb{C}^N$ and
\[
\kappa=(\kappa^+,\dots,\kappa^+_{[N/2]},\kappa^-,\dots,\kappa^-_{[N/2]})\in
(0,\infty)^{2[N/2]} .
\] The following eigenvalue problem
\begin{gather*}
T_\xi f=i\langle\lambda,\xi\rangle f,
\qquad
f(0)=1,
\qquad
\forall\, \xi\in\mathbb{C}^N,
\end{gather*}
has a~unique analytic $M_{\kappa}(\lambda,x)$ given by
\begin{gather*}
M_{\kappa}(\lambda,x)=e^{i\langle\lambda,h(0,x)\rangle}\prod_{j=1}^{[N/2]}M_{\kappa^-_j}\left(\frac{i}
{2}(\lambda_{2j-1}-\lambda_{2j})(x_{2j-1}-x_{2j})\right)
\\
\phantom{M_{\kappa}(\lambda,x)=}{}
\times
M_{\kappa^+_j}\left(\frac{i}{2}(\lambda_{2j-1}+\lambda_{2j})(x_{2j-1}+x_{2j})\right).
\end{gather*}
\end{proposition}

\subsection*{Acknowledgements}

This research is supported by NPST Program of King Saud University, project number 10-MAT1293-02.
I~would like to thank the editor and the anonymous referees for their helpful comments and remarks.

\pdfbookmark[1]{References}{ref}
\LastPageEnding

\end{document}